\documentclass[12pt]{article}

\usepackage[utf8]{inputenc}   
\usepackage[T1]{fontenc}      
\usepackage[english]{babel} 
\usepackage{amsfonts} 
\usepackage{amsmath}
\usepackage{latexsym} 
\usepackage{lmodern}  
\usepackage[all]{xy}
\usepackage{color}

\SilentMatrices
\nonstopmode   
\DeclareMathOperator{\Sub}{{\bf Sub}}
\DeclareMathOperator{\un}{{\bf 1}}
\newcommand \EE {{\cal E}}
\newcommand \JJ {{\cal J}}
\newcommand \PP {{\cal P}}
\newcommand \TT {{\cal T}}
\newcommand \N {\mathbb N}
\newcommand \rar {\rightarrow}
\newcommand \ci {\circ}

\newcommand \adj {\dashv}
\newcommand \separ[1] {\mbox{\hspace{4mm}\mbox{#1}\hspace{4mm}}}
\newcommand \ba[1] {\overline{#1}}
\newcommand \et {\wedge}
\renewcommand \implies {\Rightarrow}
\newcommand \ou {\vee}
\newcommand \equivalent {\Leftrightarrow}
\newcommand \proves {\,\vdash}
\newcommand \interp[1] {{\lfloor}#1{\rfloor}}
\newcommand \zpE {{\zeta\proves E}}
\definecolor{grey}{rgb}{0.5,0.5,0.5}
\newcommand \greyed[1] {{\color{grey}#1}}
\newcommand \gJJ[2] {\greyed{\JJ_{#1}(}#2\greyed{)}}
\newcommand \JxX[1] {\gJJ{x\in X}{#1}}
\newcommand \JzE[1] {\gJJ{\zpE}{#1}}

\newcommand \ignore[1] {}
\newenvironment{liste}{\begin{list}{$\bullet$}{\itemsep1mm\topsep0mm\partopsep0mm}}{\end{list}}

\begin{document}

\begin{center}
\bf\Large On~dependent~conjunction~and~implication
\end{center}
\vspace{1cm}
\begin{center}
\begin{tabular}{c}
{Matthieu Herrmann and Alain Prouté}\\
\end{tabular}
\end{center}

   \begin{abstract}
   In the context of our work on a new proof assistant, we were led
   led give a theoretical model of conjunctions $E\et F$ and implications $E\implies F$ 
   where $F$ is meaningful only when $E$ is true, a situation that is very often encountered
   in everyday mathematics, and which was already considered by several type theorists. 
   We present a version of these concepts that is more oriented towards usual mathematics
   than towards type theories, by using an extension of Lawvere's definition
   of quantifiers. Also, our presentation stresses the importance of proof-irrelevance 
   and of the principle of description in this phenomenon.
   We explain why this dependency is obtained through the use of a ``hidden'' variable, and more generally the
   links of these concepts with the vernacular language of mathematics, which is actually our main motivation. 
   In fact, it appears that this apparently anodine subject leads to consider several deep questions in
   the foundation of mathematics. 
   Despite its links with topos theory, this article is readable by non-specialists. 
   \end{abstract}

{\small   
\begin{flushright}
\end{flushright}
}
   
\section*{Introduction}
Conjunction, also known as the ``logical and'', and widely denoted $\et$, is generally defined
as a commutative operation. In other words, 
for any statements $E$ and $F$, it is assumed (or proved) that $E\et F$ is equivalent to $F\et E$. 

However, there are very elementary and simple examples showing that this is not always the case. 
For example, let $A$ be a subset of the set $\N$ of the natural integers, and consider the statement
$A\not=\emptyset\et\inf(A)=0$, where $\inf(A)$ denotes the infimum (greatest lower bound) of $A$. 
This statement is clearly meaningful, should it be true or not. On the contrary, the 
statement $\inf(A)=0\et A\not=\emptyset$ is meaningless, because $\inf(A)$ is not well defined if we do not
know that $A$ is nonempty. Notice that $A\not=\emptyset\implies\inf(A)=0$ is meaningful and 
$\inf(A)=0\implies A\not=\emptyset$ meaningless for the same reason. 

It is hard to find a treatement of this problem in the literature. 
It seems that the first occurrence of a theoretic analysis of this fact can be found in
a series of lectures given by P. Martin-Löf in 1980 (see \cite{MartinLof84}), within the setting of 
intuitionistic type theory.
The idea is mainly based on the fact that statements can be identified with sets (or types) in
view of the analogy between the two situations~: ``the element $a$ is of type $X$'' and ``the expression $p$ is
a proof of the statement $E$'', which is now widely popularized under the name of the ``Curry-Howard correspondance'', to
which Martin-Löf actually also refers. Subsequent works have appeared modelling these ``dependent logical connectives'', 
such as T. Coquand and
G. Huet in 1986 \cite{CoquandHuet86} (see also M. Hyland and A. M. Pitts \cite{HylandPitts87}), or D. Pavlovi\'c
\cite{Pavlovic} in 1990, among others.

We turn now to a vocabulary that is more familiar to mathematicians. 
We propose an explanation of dependent conjunction and implication,
which is not only a modelling of them, but also shows how much their existence is natural. It is indeed exactly as
natural as the existence of the quantifiers $\exists$ and $\forall$, simply because dependent conjunction and
implication play the same role with respect to assumptions as quantifiers with respect to declarations. It seems
to us that the best way to understand this is to use one of the most beautiful discoveries
of William Lawvere, namely his definition of quantifiers by way of adjoint increasing maps.(\footnote{A pair of 
    adjoint \emph{decreasing} maps (almost the same concept) is also known as a ``Galois correspondence''.})  

It appears that a dependent conjunction $E\et F$ (and similarly for an implication) 
should be denoted $(\zpE)\et F$ if we want to be completely
explicit, where $\zpE$ is a kind of declaration of $\zeta$ (actually a named assumption) and where $F$ contains free
occurrences of the variable $\zeta$. Of course, in usual mathematics, this variable $\zeta$ is not visible, and we also explain
on the one hand how $F$ can have such a free occurrence of $\zeta$, and on the other hand why $\zeta$ remains invisible
(hidden) in $F$, and consequently not explicitly declared before $E$.  

The content of this article holds for both intuitionistic and classical mathematics, and is part of
type theory and topos theory's more or less explicit folklore. We hope to add some 
originality in the way we have chosen to present it, and to make the subject accessible to a wider audience. 

The authors want to thank John Newsome Crossley for his careful reading, his suggestions and many corrections in our
use of the English language.

{\parskip-3mm\tableofcontents}

\section{Prerequisites}
The meaning of logical connectives was historically defined in two different ways~: the ``logical'' one,
which defines under what conditions a given statement is true (or false), possibly using truth tables (in other words, 
    how to compute a truth value for a statement), 
and the ``operational'' one, which defines the 
meaning of logical connectives in terms of proofs. The second manner was initiated around 1930 by the intuitionists
(and also holds for classical logic(\footnote{This is just because all intuitionistic proof principles are also valid
    in classical logic.})) through the so-called Brouwer-Heyting-Kolmogorov interpretation, 
and also by G. Gentzen in his work on natural deduction and the sequent calculus.
This operational point of view was more recently very elegantly reformulated using the categorical concept 
of adjoint functors (in the present case, adjoint increasing maps), in 
particular within the setting of topos theory.

For example, (ordinary) conjunction can be defined as the right adjoint to the duplication (diagonal)
operation.  
In order to see this, consider two preordered sets $X$ and $Y$,(\footnote{In other words, sets with a binary relation $\leq$ that is 
    reflexive and transitive.}) and two maps $f:X\rar Y$ and $g:Y\rar X$. We
say that ``$f$ is left adjoint to $g$'' (or equivalently that ``$g$ is right adjoint to $f$'') if~:
$$
\forall x\in X\ \forall y\in Y\ f(x)\leq y \equivalent x\leq g(y)
$$
This fact is denoted $f\adj g$.(\footnote{Warning~: the symbol $\adj$ is not a reversed notation for $\proves$. The
    two symbols have completely distinct uses (but they sometimes appear together in a single expression).}) 
This definition has several elementary consequences. First of all, the inequality
$x\leq g(f(x))$ is true for any $x$ in $X$ (just replace $y$ by $f(x)$ in the above statement, and use the reflexivity
    of the preorder relation). Symmetrically, we have $f(g(y))\leq y$ for any $y$ in $Y$. These two inequalities
    are known as the ``unit'' and the ``co-unit'' of the adjunction $f\adj g$. From this and the transitivity
    of preorder relations we can deduce that $f$ and $g$ are increasing functions, 
    and that if we have $f\adj g$ and $f\adj h$, then $g$ and $h$ are ``equivalent'' in the sense that
for any $y$ in $Y$, we have $g(y)\leq h(y)$ and $h(y)\leq g(y)$ (so that they are actually equal if the preorder on $X$
    is an order). Of course we also have the symmetric fact that $f\adj h$ and $g\adj h$ entail that $f$ and $g$
are equivalent (in the same sense).  
Furthermore, if $f\adj h$ and $g\adj k$, and if $g\ci f$ is meaningful, then $h\ci k$ is also meaningful and we have 
$g\ci f\adj h\ci k$. 
Finally, we shall also need the fact that left adjoints commute (up to equivalence) with least upper bounds, and that right
adjoints commute (up to equivalence) with greatest lower bounds.(\footnote{In other words, 
    if $g:Y\rar X$ is a right adjoint and if $\inf_{i\in I}(y_i)$ exists, then $\inf_{i\in I}(g(y_i))$ exists and is 
    equivalent to $g(\inf_{i\in I}(y_i))$.}) 
We leave these facts as easy exercises for the reader.

Hence, if a map has an adjoint (either on the left or on the right), this adjoint is
essentially (i.e. up to equivalence) unique. As a consequence, if a map is defined 
as a left or right adjoint to a given map, it is well defined (up to equivalence).(\footnote{This does not
    imply that the adjoint exists. It is only ``unambiguously defined (up to equivalence)'' .}) This kind of 
definition is a particular case of a ``universal problem''.

Now, consider the meta-set $\EE$ of all mathematical (closed) statements. This meta-set is preordered by the relation of
deductibility, that we denote by $\leq$ in this article. In other words, $E\leq F$ means that $F$
can be deduced from $E$ (i.e. proved under hypothesis $E$). This is clearly a preorder.(\footnote{But not an order, 
    because we are very careful 
    not to confuse statements (that are ``signifiers'', a syntactic notion) with truth values 
    (the corresponding ``signified'', a semantic notion).})
For any preordered set $X$, the product set $X\times X$ is also preordered by the relation 
defined by the condition that $(x,y)\leq (u,v)$ if and only if $x\leq u$ and $y\leq v$, and we have a ``diagonal''
(or ``duplication'') map $\Delta:X\rar X\times X$ defined by $\Delta(x)=(x,x)$ (which is obviously increasing).

Now, we can define a map $\et:\EE\times\EE\rar \EE$ just by postulating that $\Delta \adj \et$. The map $\et$ is well
defined (up to provable equivalence of statements). It is easy to recognize that this map cannot be anything other 
than our usual conjunction. This is just a consequence of the fact that the usual conjunction indeed has the property
$\Delta\adj \et$, since we have~:
$$
E\leq F\et G \separ{iff}  E\leq F \mbox{ and } E\leq G  \separ{iff} \Delta(E) = (E,E) \leq (F,G) 
$$
This definition (which is well known in topos theory) of conjunction is one of the last known avatars of the 
Brouwer-Heyting-Kolmogorov interpretation (concerning conjunction). 
It is undoubtebly the most elegant thing we can do for defining (ordinary) conjunction. 

Why is such a definition ``operational''~? This is just because all the proof rules concerning conjunction
can be deduced from it, which is not the case if conjunction is defined by the ``logical'' method, i.e. in terms 
of truth values. Indeed, the adjunction $\Delta\adj\et$ means that $E\leq F\et G$ if and only if 
$E\leq F \mbox{ and } E\leq G$, in other words that it is equivalent to prove $F\et G$ under the hypothesis $E$, or to
prove separately $F$ and $G$ under the same hypothesis. This is of course a rule that anyone uses everyday when a
conjunction must be proved. Now, consider the co-unit of the adjunction, which gives~:
$$
(E\et F,E\et F) = \Delta(E\et F) \leq (E,F)
$$
and which gives the rules $E\et F\leq E$ and $E\et F\leq F$ (which are generally taken as axioms in most
    presentations of the properties of conjunction).  The reader is urged to check that similarly 
we have $\ou\adj \Delta$, and that the principle of 
    ``reasoning by disjunction of cases'' is actually a consequence of this definition.(\footnote{It should be remarqued 
that there is no dependent disjunction in mathematics. Indeed, a disjunction $E\ou F$ should be true if $E$ is false 
and $F$ true. But in this case, $F$ cannot depend on the truth of $E$.})

Actually, $E\et F$ is nothing other than the greatest lower bound of $E$ and $F$, as is the case of any concept defined
as the right adjoint of a diagonal map (as for example gcd using the divisibility relation, or intersection using
    the inclusion of subsets relation). 
However, the notion of adjoint maps is much more general than the notion of
greatest lower bound. Of course, commutativity of the (ordinary) conjunction is a consequence of the commutativity
of the diagonal~$\Delta$. 

Because we shall use it in another section below, we also recall the definition of the usual implication via
adjunctions. Given a statement $F$, we consider the map $E\mapsto E\et F$ from $\EE$ to $\EE$. By definition, the map
$G\mapsto F\implies G$ is a right adjoint to the previous one. In other words, we have~:
$$
E\et F \leq G \separ{iff} E \leq F\implies G
$$
This means that in order to prove the implication $F\implies G$ under the hypothesis $E$ it is enough to prove $G$ under
the hypothesis $E\et F$. This is of course what everyone does when an implication must be proved, and is known as the
``auxiliary hypothesis method''.(\footnote{\label{fnote1}Purists will argue that the method consists in considering $F$ as an \emph{extra 
    hypothesis}, not to replace the hypothesis $E$ by the hypothesis $E\et F$. They are right. Passing from one point of view to the
    other one can of course be formalized, and we will need to do it in section \ref{sec:depI}. For the time being, we need to present these
 concepts as we do, because otherwise adjunctions are not as simply applicable.})

The co-unit of this adjunction gives~:
$$
(E\implies F)\et E \leq F
$$
which is the principle known as ``modus ponens''. 

Another advantage of this method, is that it gives a clear and exciting definition of the quantifiers, 
as was revealed by Lawvere.

\section{Lawvere's definition of the quantifiers}
In order to express Lawvere's definition of the quantifiers, we must first discuss the notion of ``context''. In
mathematical texts, it is often the case that we ``declare'' variables. For example, we can say ``Let $x$ be a real
number.''. Such a sentence is called a ``declaration''. In this article, a ``generic'' 
declaration will be denoted $(x\in X)$, where $x$ is a symbol and where 
$X$ is a set. 

Within any part of a mathematical text, we are ``working in a context'', which
is just the consequence of the declarations made so far.(\footnote{And the ``scope'' of which we have not yet left.}) 
Let $\Gamma$ denote an arbitrary context, and let $\Gamma(x\in X)$ denote the 
    context obtained by declaring $x$ as an element of $X$ in the context $\Gamma$. In other words, $\Gamma(x\in X)$ is the context obtained by
    ``enriching'' $\Gamma$ by a new declaration.(\footnote{Formally, contexts are generated by the two rules~: there is an 
        ``empty'' context, and for any context $\Gamma$, $\Gamma(x\in X)$ is a context.}) 
    We require that a variable cannot be declared twice in the same context. In other words, the
    context $\Gamma(x\in X)$ is invalid if $x$ is already declared in $\Gamma$.  
    
    We also denote by $\EE_\Gamma$ the meta-set of all statements
    that are meaningful in the context $\Gamma$, in other words, the free variables of which are declared in $\Gamma$. 
    This meta-set $\EE_\Gamma$ is still preordered by the deductibility relation, that we denote $\leq_\Gamma$. 
    There is a canonical inclusion $\JJ_{x\in X}:\EE_\Gamma\rar \EE_{\Gamma(x\in X)}$, since if a free
    variable in a statement $E$ is declared in $\Gamma$, it is {\em a fortiori} declared in 
    $\Gamma(x\in X)$. Notice the importance of the fact that a variable cannot be declared twice in the same context.
    Indeed, if a variable $x$ is already declared as an element of some set $Y$ in $\Gamma$, it can be the case that a
    statement $E$ which is meaningful in $\Gamma$ becomes meaningless in $\Gamma(x\in X)$. This phenomenon is analogous
    to the well known ``variable capture''.   
    
    Since $\JJ_{x\in X}$ is a canonical
    inclusion, any statement $E$ is identical to $\JJ_{x\in X}(E)$ from a syntactic point of view. However, it is in many
    places very important to still see this operator $\JJ_{x\in X}$. This is why, instead of writing $\JJ_{x\in X}(E)$ as $E$, 
    we write it in grey~: $\JxX{E}$. 
    
    Now, if $a$ is \emph{any} expression
    representing an element of $X$ in the context $\Gamma$, and if $E\in\EE_{\Gamma(x\in X)}$, the process of replacing
    all (free) occurrences of $x$ in $E$ by the expression $a$ produces a statement interpretable (i.e. meaningful) in the context $\Gamma$, 
    in other words, an
    element of $\EE_\Gamma$, that we denote by $E[a/x]$ (read ``$E$ where $a$ replaces $x$''). Hence, we have a map 
    $[a/x]:\EE_{\Gamma(x\in X)}\rar \EE_\Gamma$ (for each such $a$), and it is clear that this map is a retraction for $\JJ_{x\in X}$ 
    (i.e. that $[a/x]\ci\JJ_{x\in X} = 1_{\EE_\Gamma}$). Furthermore, this map is increasing (in other words, 
        if $E\leq_{\Gamma(x\in X)} F$, then $E[a/x]\leq_\Gamma F[a/x]$)
    because the replacement can also be performed within proofs.(\footnote{We do not give any precise definition
        of the replacement. We assume only that it is a retraction for $\JJ_{x\in X}$ and that it is increasing. 
        These two properties are of course easy consequences of the ``usual'' definition of the replacement, and of the definition
        of proofs (that we do not explicitly need here). 
        
        \emph{Important warning}~: The inclusion 
        $\JJ_{x\in X}:\EE_\Gamma\rar\EE_{\Gamma(x\in T)}$ allows to identify $\EE_\Gamma$ to a subset of $\EE_{\Gamma(x\in X)}$, but
        \emph{not as a preordered subset}. Indeed, we can have $\JxX{E}\leq_{\Gamma(x\in X)}\JxX{F}$, and not have $E\leq_\Gamma F$, 
        because the fact of declaring $x$ in $X$ can allow to prove things which are not provable in the context $\Gamma$. 
        Think for example of what happens if you declare an element in the empty set.})
    
Lawvere defines the quantifiers by postulating that~:(\footnote{Maybe the word ``defines'' needs an explanation. The reality
    is that such adjunctions can be used to define the deductibility relation itself, and this provides a meaning for the logical
    connectives as a by-product. See \cite{cours} Chapitre 1, for a detailed explanation of this fact. Lawvere did 
    not actually give the definition
    exactly in this form. His ``definition'' (or more accurately, his characterization) looks much more in the form of a pair of internal 
    (in the sens of the ``internal logic'' of topos theory) 
    adjoints to the ``inverse image'' arrow $\PP(Y)\rar\PP(X)$ (for a given arrow $X\rar Y$), where $\PP(X)$ and $\PP(Y)$ are 
    seen as ``internally ordered objects''. There is also an ``external'' version of this using $\Sub$ 
    (the subobject functor) instead of $\PP$ (the power object functor). 
    In some sens, Lawvere's actual definition uses subsets where we are 
    using statements, but this is of course equivalent in view of the one-to-one correspondance between subobjects 
    (or subsets) and characteristic arrows. See Lawvere~\cite{Lawvere}.})(\footnote{Writing
    the declaration after a quantifier as subparts is not usual practice. However, we want to stress the fact that the 
    declaration $x\in X$ plays the role of a \emph{parameter} for the function $E\mapsto \forall_{x\in X}\ E$.
    This is also coherent with the notation $\JJ_{x\in X}$.})
$$
\exists_{x\in X}\adj \JJ_{x\in X}\adj \forall_{x\in X}
$$
Let us check that our ``usual'' quantifiers have this property, so that they are actually equivalent to Lawvere's
quantifiers. In the case of the universal quantifier, Lawvere's definition writes~:
$$
\JxX{E}\leq_{\Gamma(x\in X)} F  \separ{iff}  E\leq_\Gamma \forall_{x\in X}\ F
$$
What this says is that in order to prove $\forall_{x\in X}\ F$
under the hypothesis $E$, we can just declare $x$ in $X$ and prove $F$ under this same hypothesis $E$ (and that 
    the converse is also true). This is clearly
what we do everyday. 

The co-unit of the adjunction writes~:
$$
\JxX{\forall_{x\in X}\ F} \leq_{\Gamma(x\in X)} F
$$
which is the most general case of particularization of a universally quantified statement. It says that
if $\forall_{x\in X}\ F$ is true, then $F$ is true, provided of course that we interpert this in the context $\Gamma(x\in X)$ (otherwise, 
    $F$ would anyway be meaningless). Applying the replacement map $[a/x]$ to both sides, we get~:
$$
\forall_{x\in X}\ F \leq_\Gamma F[a/x]
$$
which is the ``usual'' particularization principle. 

A symmetric analysis can be done for the existential quantifier. 
Precisely, the adjunction writes~:
$$
\exists_{x\in X}\ E \leq_{\Gamma} F \separ{iff} E\leq_{\Gamma(x\in X)} \JxX{F}
$$
This shows how to use an existence hypothesis. Indeed, in order to prove $F$ under the hypothesis $\exists_{x\in X}\ E$,
we can first declare $x\in X$, and then prove $F$ under the hypothesis~$E$. This is again what everybody does
instinctively. Applying the replacement $[a/x]$ to the unit of the adjunction yields~:
$$
E[a/x] \leq_\Gamma \exists_{x\in X}\ E
$$
which is the usual ``exhibition principle'' used for proving an existence.

\section{Dependent conjunction and implication}\label{sec:depI}
We now arrive at the heart of our subject. 
Contexts are constructed not only via declarations, but also via ``assumptions''. 
For example, we
can find the following sentence within a mathematical text~: ``Let $x$ be a real number, and assume that $x>0$.'', i.e. a
declaration followed by an assumption. Obviously, both are modifying the context in the intuitive sense, 
    simply because they both provide new means for proving. 

Let's use the notation $\zpE$ (read ``$\zeta$ proves $E$\,'') for a generic assumption, where $E$ is a
statement, and $\zeta$ the ``name'' of the assumption.
 As previously with declarations, we consider contexts of the form $\Gamma(\zpE)$ obtained by enriching $\Gamma$ by
 the ``assumption'' $\zpE$,(\footnote{Declarations and assumptions can be mixed in a single context.}) 
 and we have a canonical inclusion $\JJ_{\zpE}:\EE_\Gamma\rar\EE_{\Gamma(\zpE)}$.(\footnote{In general, 
     assumptions are anonymous~: ``Assume $E$.'', 
    but can also be named such as in~: ``Assume $E$ ($\zeta$).'', where $\zeta$ is a symbol (a name, but in practice, 
        a number is generally used) by which we can later refer to
    this assumption. The reason why assumptions are generally anonymous is a consequence of a fundamental principle
    of mathematics, called ``proof-irrelevance'',  
    which is discussed in section \ref{sec:uniqueness}. For our explanations, 
    we need to be explicit, hence the notation 
    $\zpE$. Furthermore, if $p$ is any proof of $E$ in the context $\Gamma$, the replacement operation 
$[p/\zeta]:\EE_{\Gamma(\zpE)}\rar\EE_\Gamma$ is a retraction for $\JJ_{\zpE}$ and is increasing. In this article, 
    we do not give any precise definition of the notion of proof. This is not necessary for our purpose, and we will 
    actually not use this kind of replacement.})

Imitating Lawvere, we introduce, for any statement $E$, two ``declarative operators'', 
$$F\mapsto (\zpE)\et F \separ{and} F\mapsto (\zpE)\implies F$$
which are actually two maps from $\EE_{\Gamma(\zpE)}$ to $\EE_\Gamma$, by postulating~:
$$
(\zpE)\,\et \adj \JJ_{\zpE} \adj (\zpE)\implies
$$
Before proving anything from this definition, we must capture the fact that $E$ is true in the context 
$\Gamma(\zeta\proves E)$ (and even a little more than this). 
Nothing up to here can entail
this fact, because there is \emph{a priori} no link between the deductibility relation and the concept represented by the
sign $\proves$, which actually has no precise meaning up to now, and is mainly just a syntactic gadget. In other words,
we must state a principle that establishes the link between ``deductibility'' and ``assumption''. We propose the following
\emph{special rule}~:
$$
E\et F \leq_\Gamma G \separ{iff} \JzE{F}\leq_{\Gamma(\zpE)}\JzE{G}
$$
(where $\et$ is ordinary conjunction) 
which says that it is the same thing to deduce $G$ from $E\et F$ or to deduce $G$ from $F$ alone after having assumed
$E$. It is actually the formalisation of the difference refered to in footnote \ref{fnote1}.(\footnote{A similar rule 
    exists in Gentzen's work (left $\et$-rule), which establishes the link between the conjunction and the coma in the left members of
    sequents. This coma may be considered as an external version of the conjunction, and indeed, concatanation of
    assumptions in a context is also a kind of ``external'' conjunction. From this special rule, we can derive
the rule~:
$
H\leq_{\Gamma(\zpE)} \JzE{E}
$
for any statement $E$ meaningful in $\Gamma$ and any statement $H$ meaningful in $\Gamma(\zpE)$. Indeed, it is enough to 
replace $F$ by $\top$ (a greatest element in $\EE_\Gamma$), and $G$ by $E$ in the special rule, and to use the 
facts that $E\et\top \simeq E$ and
that 
$H\leq_{\Gamma(\zpE)}\top\simeq \JzE{\top}$, because indeed, $\JJ_\zpE$, as a right adjoint, 
    preserves greatest elements which are the greatest lower bounds of the empty subset. This derived rule just says 
    that if you assume $E$, then $E$ is true.})

Now, let us examine the consequences of the above definition. 
Concerning the right adjoint (the ``dependent implication'') we have the equivalence~:
$$
\JzE{F}\leq_{\Gamma(\zpE)} G \separ{iff} F \leq_\Gamma (\zpE)\implies G
$$
This means that proving the implication
$(\zpE)\implies G$ under the hypothesis $F$ is the same as proving $G$ under the same hypothesis $F$, but
after having assumed $E$. This is clearly the same thing as the usual ``auxiliary hypothesis method'',(\footnote{And here, 
    we are much closer to the usual meaning of this method than in the introduction. See footnote~\ref{fnote1}.}) and this just
shows that this method is compatible with the dependency of $G$ on the truth of $E$. 

Notice that the co-unit of this adjunction says~:
$$
\JzE{(\zpE)\implies G} \leq_{\Gamma(\zpE)} G
$$
which means that if you have the hypothesis $(\zpE)\implies G$ in 
a context where $E$ is assumed, you can deduce $G$. This is the most general dependent version of ``modus ponens''. 
There is no hope to have exactly the form  $(E\implies G)\et E \leq G$ because the left hand side should be interpretable in a 
context not declaring $\zeta$, where
$G$ is not meaningful (because it depends on~$\zeta$).(\footnote{However, if we have a proof $p$ of $E$ in the context $\Gamma$, 
    applying $[p/\zeta]$ to both sides, we get $(\zpE)\implies G\leq_\Gamma G[p/\zeta]$, which looks like another variation on modus ponens.})
However, we have the following inequality~:(\footnote{This inequality was suggested to us by Paul-André Melliès during a
    talk on this subject. Actually, this is a particular case of the fact that if $F$, $G$ and $H$ are functors such that 
    $F\adj G\adj H$, then we have the natural transformation $\varepsilon\ci G\varepsilon H$ from $GFGH$ to the identity 
    functor (where the two $\varepsilon$ are the co-units of the two adjunctions).})
$$
\JzE{(\zpE)\et \JzE{(\zpE)\implies G}}\leq_{\Gamma(\zpE)} G
$$
Indeed, we first have $(\zpE)\et \JzE{(\zpE)\implies G}\leq_{\Gamma} (\zpE)\implies G$, by the co-unit of the adjunction
$(\zpE)\,\et \adj\JJ_{\zpE}$, and since $\JJ_{\zpE}$ is increasing, we have~:
$$
\JzE{(\zpE)\et \JzE{(\zpE)\implies G}}\leq_{\Gamma(\zpE)} \JzE{(\zpE)\implies G}
$$
Now, we also have $\JzE{(\zpE)\implies G}\leq_{\Gamma(\zpE)} G$, by the co-unit of the other adjunction. We shall show
in section \ref{sec:comp} that $(\zpE)\et \JzE{(\zpE)\implies G}$ is equivalent to $E\et ((\zpE)\implies G)$ (where $\et$ is ordinary
    conjunction), so that we also have the simplified form~:
$$
\JzE{E\et ((\zpE)\implies G)}\leq_{\Gamma(\zpE)} G
$$
which is probably the formulation that is the closest possible to ordinary modus ponens.

Concerning the left adjoint (``dependent conjunction'') we have the
equivalence~:
$$
(\zpE)\et F \leq_\Gamma G\separ{iff} F\leq_{\Gamma(\zpE)} \JzE{G}
$$
which means that in order to use the hypothesis $(\zpE)\et F$ in order to prove $G$, 
we can just use only the hypothesis $F$, but after having assumed $E$. This is again indeed what
we do daily in a natural and instinctive way.(\footnote{Notice the similarity between this equivalence and the special 
    rule. In some sense, the special rule says that the ordinary conjunction also has this property.})

The unit of this adjunction says~:
$$
F \leq_{\Gamma(\zpE)} \JzE{(\zpE)\et F}
$$
which means that $(\zpE)\et F$ can be deduced from $F$, provided that $E$ is assumed, which is 
again an everyday proof method (actually, the \emph{canonical} way for proving a dependent conjunction).

Notice that the non dependent counterpart of this last fact is just 
the fact that in order to prove $E\et F$ under some hypothesis, we can prove $E$ and then $F$ under 
the same hypothesis. In the dependent situation this works similarly except that $E$ must be proved first, 
because the fact that $E$ is true is necessary to make $F$ 
meaningful.

\section{Comparison with ordinary conjunction and implication}\label{sec:comp}
A natural question is that of the equivalence of $(\zpE)\et \JzE{F}$ and ordinary conjunction $E\et F$, and 
similarly for implication. Here,
    $F$ is meaningful in the same context as $E$, so that $\JzE{F}$ does 
    not ``actually depend'' on the truth of $E$.(\footnote{Be careful not to confuse the notion of ``non-dependent'' which 
applies to ordinary conjunction and implication, and the notion of ``not actually dependent'', which applies to
dependent conjunctions
and implications of the forms $(\zpE)\et \JzE{F}$ and $(\zpE)\implies \JzE{F}$.})  

In order to avoid any confusion, we recall that a dependent conjunction is always denoted $(\zpE)\et F$ even if $F$ does
not actually depend on the truth of $E$ (i.e., even if $F$ is in the image of $\JJ_{\zpE}$). On the contrary, of course,
ordinary conjunction is denoted $E\et F$. The same rule applies to the two kinds of implication.

For conjunction, we have successively~:
$$
\begin{array}{rlll}
E\et F           &\leq_\Gamma        & E\et F            &\mbox{reflexivity}\\
\JzE{F}          &\leq_{\Gamma(\zpE)}& \JzE{E\et F}      &\mbox{(special rule)}\\
(\zpE)\et\JzE{F} &\leq_\Gamma        & E\et F            &\mbox{($(\zpE)\,\et\adj\JJ_\zpE$)}
\end{array}
$$
For the converse inequality, we have~:
$$
\begin{array}{rlll}
\JzE{F}        &\leq_{\Gamma(\zpE)} & \JzE{(\zpE)\et\JzE{F})}    &\mbox{(unit of the adjunction)}\\
E\et F         &\leq_\Gamma         & (\zpE)\et\JzE{F}           &\mbox{(special rule)}
\end{array}
$$

For implication, we have successively~:
$$
\begin{array}{rlll}
(E\implies F)\et E      &\leq_\Gamma        & F              &\mbox{(modus ponens)}\\
\JzE{E\implies F}       &\leq_{\Gamma(\zpE)}& \JzE{F}        &\mbox{(special rule)}\\
E\implies F             &\leq_\Gamma& (\zpE)\implies \JzE{F} &\mbox{($\JJ_\zpE\adj (\zpE)\implies$)}     
\end{array}
$$
and for the converse inequality~:    
$$
\begin{array}{rlll}
\JzE{(\zpE)\implies \JzE{F}}     &\leq_{\Gamma(\zpE)}& \JzE{F}        &\mbox{(co-unit of the adjunction)}\\
((\zpE)\implies \JzE{F}) \et E   &\leq_\Gamma        & F              &\mbox{(special rule)}\\
(\zpE)\implies \JzE{F}           &\leq_\Gamma        & E\implies F    &\mbox{($A\mapsto A\et E\adj F\mapsto E\implies F$)}
\end{array}
$$
Hence, the dependent conjunction and implication are equivalent to ordinary conjunction and implication when their
second operand does not ``actually depend'' on the truth of the first one. This is of course an important and reassuring
fact. 

Possibly surprising may be the fact that usual conjunction comes as a right adjoint (and consequently 
    is a ``multiplicative'' connective), 
whereas dependent conjunction comes as a left adjoint (and hence is an ``additive'' connective). In other words, 
the meaning of the usual conjunction is defined in terms of its behaviour when it is in the position 
of the conclusion, whereas the
meaning of dependent conjunction is defined in terms of its behaviour when it is in the position of the 
hypothesis, as we saw above. A fact which
can demystify this point for the reader is the similar fact that an indexed disjoint union of sets, such as
$$
\coprod_{x\in X}\ Y
$$
(a notion of an additive nature, which can actually be defined through the use of a left adjoint) 
can be identified with the cartesian product $X\times Y$ (a notion of a multiplicative
    nature) when $Y$ does not actually depend on $x$, which is after all just a more elaborate example than 
$2+2+2=3\times 2$. 

However, we can temper the mysterious aspect of this by remarking that ordinary conjunction as a right adjoint
is a binary operation, whereas dependant conjunction is in no way a binary operation. It looks much like a
declarative operation, i.e. it constructs expressions containing a declaration and a ``body'' which is the scope of this
declaration. There are many such declarative operations in mathematics, for example~:
$$
\forall_{x\in X}\ E \hspace{6mm}
\{x\in X \ |\ E\} \hspace{6mm}
(x\in X) \mapsto E  \hspace{6mm}
\coprod_{x\in X}\ E  \hspace{6mm}
\mbox{etc}\dots
$$
(where $E$ is the body). 
So, the dependent conjunction is a quite different operation from ordinary conjunction.  
In a conclusion, we may be troubled by the fact that the two conjunctions have the same name
(which is nevertheless justified by what we proved above). 

The fact that conjunction should preferably be considered as a left adjoint (or an additive connective) rather than 
a right adjoint (or a multiplicative connective) 
could have been (but actually was not) suggested in the first half of the twentieth century. Indeed, the 
Brouwer-Heyting-Kolmogorov interpretation, which is a definition of the meaning 
of the logical connectives in terms of proofs, says that~:(\footnote{For a discussion of this interpertation stressing the 
    fact that it hides a confusion between the signifier and the signified, see \cite{cours}, pages 25-26. This confusion 
    disappears if one replaces ``proof'' by ``warrantor''; see section \ref{sec:uniqueness}.})
\begin{liste}
\item a proof of $E\et F$ is a \emph{pair} $(p,q)$ where $p$ is a proof of $E$ and $q$ a proof of $F$, 
\item a proof of $E\ou F$ is a \emph{pair} $(i,p)$ where $i$ is either $0$ or $1$, and $p$ a proof of $E$ if $i=0$ and a
proof of $F$ if $i=1$, 
\item a proof of $E\implies F$ is a \emph{method}(\footnote{By ``method'', you can understand ``algorithm''.}) for producing
a proof of $F$ from a proof of $E$, 
\item a proof of $\forall_{x\in X}\ E$ is a \emph{method} for producing a proof of $E$ from any $x$ in $X$, 
\item a proof of $\exists_{x\in X}\ E$ is a \emph{pair} $(x,p)$ where $x$ is an element of $X$, and $p$ a proof of
$E$.
\end{liste}
As one can see, additive connectives are those for which a proof is a pair, whereas multiplicative connectives are those
for which a proof is a method, \emph{provided that the conjunction is considered as an additive connective}. However, it
should be remarked that pairs which are proofs of $E\ou F$ and $\exists_{x\in X}\ E$ are clearly ``dependent pairs'', in
the sense that the statement proved by the second component of the pair depends on the \emph{value} of the first component, and
it is certain that, despite the fact that they were necessarily aware of this dependency, the early intuitionists 
always considered pairs that are proofs of a conjunction as ``ordinary pairs''. In the work of Gentzen, there is also
no mention of a dependent conjunction. The reason for this could be that for some reason, these mathematicians were not
much concerned by the principle of description (see section \ref{sec:desc}), or equivalently, but maybe more accurately, that they did not
take into account in their models of reasoning the fact that the definition of a mathematical object can depend on the truth of
a statement.

\section{From types to sets}
The definition given in 
section \ref{sec:depI} of dependent conjunction and implication provides a clear explanation of another phenomenon
that we now consider. 

The authors of this article were working on a compiler for 
a typed system similar to those familiar to topos theorists, such as the type theory one can find in
Lambek and Scott \cite{LS}, but with a type constructor $W$ (taking a statement as its unique operand) with the
    meaning that data of type $W(E)$ are ``warrantors'' of $E$, i.e. objects that ``warrant'' the truth 
    of~$E$.(\footnote{More on this notion in section \ref{sec:uniqueness}.})
    At the ``low level'' of this system, the declaration following a quantifier has
the form $x:T$ where $T$ is a type (not a set~!), and where the symbol $:$ plays a role similar to $\in$, but means ``is
of type'' instead of ``belongs to''. In this system, sets are defined as data whose type has the form $\PP(T)$, where
the type constructor $\PP$ is of course a syntactic version of the power object functor of topos theory. More precisely, a
datum of type $\PP(T)$ is called a ``set hosted by $T$''. We also say that $T$ is the ``host'' of~$X$.  

At the ``high level'' of the system, types are invisible and only sets are manipulated. In particular, the quantifiers
have the form $\forall_{x\in X}\ E$ and $\exists_{x\in X}\ E$ where $X$ is a set (not a type~!). The compiler of the system must
translate this high level formulation into a low level one. The formulas for this transformation are~:
$$
\begin{array}{rcl}
\forall_{x\in X}\ E &:=& \forall_{x:T}\ (x\in X\implies E)\\
\exists_{x\in X}\ E &:=& \exists_{x:T}\ (x\in X\et E)
\end{array}
$$
where $T$ is the host of $X$, and where $x\in X$ in the right hand members is not a declaration but a statement (whereas
    it is a ``high level'' declaration in the left hand members). 

Why these definitions should be written like this is intuitively obvious. However, it is not so easy to give an explanation
of why they \emph{must} be written like this. This can be done as follows.

First of all, we need a definition for the high level declaration $x\in X$. It is defined as the low level declaration
$x:T$ (where $T$ is the host of $X$), 
accompanied by the assumption $x\in X$.(\footnote{This is just because $X$ is to be considered intuitively as a ``part'' of $T$.})
In other words, the high level context enrichment by $(x\in X)$
must be translated into the low level ``two step'' context enrichment as $(x:T)(\zeta\proves x\in X)$.

Now, since the adjoint of a composition is the (reversed) composition of the adjoints, we have~:(\footnote{Recall that the 
    notations $(\zeta\proves x\in X)\,\et$ and $(\zeta\proves x\in X)\implies$ are abbreviations for
    $E\mapsto (\zeta\proves x\in X)\,\et\, E$ and $E\mapsto (\zeta\proves x\in X)\implies E$.})
$$
\exists_{x:T}\ci (\zeta\proves x\in X)\,\et
\adj
\JJ_{\zeta\proves x\in X}\ci\JJ_{x:T}
\adj
\forall_{x:T}\ci (\zeta\proves x\in X)\implies
$$
which is just the wanted explanation. 

Notice that this explanation is made possible by the fact (among other things) that dependent conjunction is defined as a
{\em left} adjoint. Also notice that in the expressions $\forall_{x:T}\ (x\in X\implies E)$ and $\exists_{x:T}\ (x\in
X\et E)$, the implication and the conjunction are obviously dependent, since $E$ can be meaningless if $x$ does not
belong to $X$.

We can apply the same method to similar situations. We give just one example. Denote by $\TT(T)_\Gamma$ the
meta-set of all terms of type $T$ meaningful in the context $\Gamma$.(\footnote{$\EE_\Gamma$ could be considered as the special
    case $\TT(\Omega)_\Gamma$, where the object $\Omega$ is the subobject classifier.}) The meta-set $\TT(\PP(T))_\Gamma$ of all expressions representing sets hosted by $T$
    is preordered by the relation of inclusion, denoted $\subset$.(\footnote{Indeed, the sets hosted by $T$ are ordered by 
        inclusion, whereas the expressions which represent them are only preordered by inclusion, since several distinct
        expressions can represent the same set. In logic, we should never confuse a signifier with the corresponding 
        signified, as we were taught by Gottlob Frege~\cite{Frege}.}) We use the more precise notation $\subset_\Gamma$ in 
    order to stress the fact that the inclusion is valid
    in the context $\Gamma$ (which entails that both operands of $\subset_\Gamma$ are meaningful in the context 
        $\Gamma$). We have the two ``low level'' maps~:
$$
\xymatrix@C=32mm{
\EE_{\Gamma(x:T)}\ar@/^/[r]^-{E\mapsto\{x:T\ |\ E\}}&\TT(\PP(T))_\Gamma\ar@/^/[l]^-{A\mapsto x\in A}
}
$$
and similarly, two ``high level'' maps by replacing $x:T$ by the declaration $x\in X$, and $\PP(T)$ by $\PP(X)$. 
Notice that 
$E\mapsto\{x:T\ |\ E\}$ is left adjoint to $A\mapsto x\in A$, since~:
$$
\{x:T\ |\ E\}\subset_\Gamma A \separ{iff} E \leq_{\Gamma(x:T)} x\in A
$$
and similarly for the high level maps. 
Now, we have~:
$$
(E\mapsto \{x:T\ |\ E\})\ci((\zeta\proves x\in X)\,\et)
\adj
\JJ_{\zeta\proves x\in X}\ci (A\mapsto x\in A)
$$
again by the composition of adjoints. Since, $\gJJ{\zeta\proves x\in X}{x\in A}$ can also be written $x\in A$, we see
that the high level comprehension $\{x\in X\ |\ E\}$ must be translated at the low level into
$\{x:T\ |\ x\in X \et E\}$.

\section{How an invisible variable can occur free in an expression}\label{sec:desc}
The declaration of $\zeta$ in $(\zpE)\et F$ and in 
$(\zpE)\implies F$\,(\footnote{$\zpE$ is of course a kind of declaration.}) would be useless if $\zeta$
had no free occurrence (should it be hidden) in $F$. Of course, we know that in usual mathematics, this $\zeta$ is
not explicitly declared and not explicitly written in the expression $F$ (which is why the phenomenon is somewhat
    mysterious). 

Returning to our example $A\not=\emptyset\et \inf(A)=0$ in the introduction, we want to write it as follows~: 
$$
(\zeta\proves A\not=\emptyset)\et \inf[\zeta](A)=0
$$
in order to make everything explicit, and in particular the fact that $\inf$ \emph{needs} $\zeta$. The purpose of this
section is to explain how a free occurrence of $\zeta$ can be present (even if invisible) in the statement $\inf(A)=0$. 

In everyday mathematics, there is a principle that we apply without being in general conscious of the fact that it 
is one of the key bricks of the system. This is the ``principle of description'', which (informally) says that if you
have proved the {\em existence} and {\em uniqueness} of a mathematical object, then this object is well defined. In topos
theory, this principle is a theorem that can be stated as follows. Given two objects $X$ and $Y$ in a topos, and an
arrow $\varphi:X\times Y\rar \Omega$,(\footnote{Recall that $\Omega$ is the subobject classifier, which plays the role
    of an ``object of truth values'', similar to the usual set of booleans $\{\mbox{true},\mbox{false}\}$, but possibly much more 
    complicated than the booleans (and which is not necessarily a set).}) such that the statement of the internal language $\forall_{x\in X}\ \exists!_{y\in
Y}\ \varphi(x,y)$(\footnote{Where as usual, $\exists!$ means ``there exists a unique''.}) is true (in the empty context), then there is one and only one arrow $f:X\rar Y$ such that the
statement of the internal language $\forall_{x\in X}\ \varphi(x,f(x))$ is true.(\footnote{It is easier to see this theorem
    as a version of the informal explanation above if you consider that $X$ represents a context (instead of the domain
        of a function), or alternatively by putting $X=\un$.}) 
The importance of the role played in mathematics by this principle is discussed in~\cite{Pr2001}. 

If we want to make everything explicit (so that mathematics look like a programming language), we introduce a so-called
``description operator'', taking as unique operand a proof(\footnote{Or more accurately a warrantor; see section \ref{sec:uniqueness}.}) 
$p$ of a statement of the 
form $\exists!_{x\in X}\ E$, 
and producing an element of $X$ (satisfying the statement $E$). This operator exists in every topos.(\footnote{It is
    essentially the arrow denoted $\sharp$ in \cite{Pr2001}.}) Notice
that considering a similar operator valid for any statement of the form $\exists_{x\in X}\ E$ (i.e. without the
    uniqueness
    condition) would be equivalent to introduce something like a ``Hilbert choice operator'' (and this one does not
exist in every topos(\footnote{A topos satisfies the ``external axiom of choice'' if all its objects are projective. This 
    is stronger than the ``internal axiom of choice'' which says that the axiom of choice is true as a statement of the internal 
    language. A Hilbert choice operator is even stronger than the external axiom of choice because it amounts to choose
    a section for every epimorphism. Actually, the original Hilbert choice operator $\varepsilon$ applied to statements and was such that 
    $\varepsilon(\exists_{x\in X}\ E)$ was well defined even if there is no $x$ such that $E$. This allowed Hilbert
    to define the meaning of $\exists_{x\in X}\ E$ as that of $E[\varepsilon(\exists_{x\in X}\ E)/x]$. This is why we
    say ``something like the Hilbert choice operator''. For the original Hilbert's work, see \cite{Heijenoort}})). 

Now, it is clear that the definition of the greatest lower bound of a nonempty subset of $\N$ looks like this~:
``{\em the unique integer such that \dots}'', in other word, it uses the principle of description, or, if we want
everything to be explicit, the description operator. This is how the variable $\zeta$ has a free occurrence in $\inf(A)$,
since the expression defining $\inf(A)$ necessarily has a subterm of the form $\delta(p)$ (where $\delta$ is the 
    description operator), where 
$p$ is a proof of a statement of the form $\exists!_{x:T}\ E$, and has at least one free occurrence of $\zeta$, since it
must use $\zeta$ as an hypothesis.

In the next section, we explain why all those things are invisible.

\section{Why a variable can be invisible}\label{sec:uniqueness}
Proof-irrelevance
has important consequences for understanding the usual language of mathematics, 
the so-called ``vernacular'' of mathematics. It explains
why assumptions are most often anonymous, why a mathematician can use a theorem without reading a proof
of it, why proofs can be written in a somewhat lazy way, in contrast with statements and terms,(\footnote{Terms are those expressions 
    representing mathematical objects in the usual sense.}) why some formal proof constructors have no counterpart in the 
everyday language of mathematics, 
and also why the variable $\zeta$ in the dependent conjunction $(\zpE)\et F$
and the dependent implication $(\zpE)\implies F$ is invisible.

Proofs are a syntactic concept, i.e. they are expressions of a
language, so that they are candidates to be ``signifiers''. As any signifier they have a corresponding signified. 
We call such a signified a ``warrantor''. But be careful, a warrantor is not just an equivalence class of proofs. The
notion of warrantor can be precisely defined in topos theory,(\footnote{A statement $E$ meaningful in a context $\Gamma$ 
    is interpreted as an arrow $\interp{E}_\Gamma$ from $\ba\Gamma$ (the object representing $\Gamma$) to $\Omega$. A
    warrantor of $E$ is just a section of the pullback of $\top:\un\rar\Omega$ along $\interp{E}_\Gamma$. This notion
    was introduced by the second author around 2007 in \cite{Luminy}.}) 
and it can be the case that a given warrantor is
represented by no proof at all.(\footnote{Which just means that in the given topos we have a true statement which is not
    (internally) provable.}) However, a warrantor of a statement $E$ always \emph{warrants} that $E$ is true (hence
    its name). Proof-irrelevance is then a consequence of the principle of ``the uniqueness of
the warrantor'', i.e. the fact that a given statement cannot have more than one warrantor (which is actually an easy consequence 
    of its topos theoretic definition).

We now use a metaphor. 
Imagine that you are alone on a desert island. Do you think that you
need a name~? For sure, you do not. You are just ``the one on the island''(\footnote{Some would argue that this 
    expression also is a ``name''. However, this kind of name could not be used for differentiating two people 
    on the island.}) 
(of course, we know which island (which plays here the role of the statement) 
    we are talking about). Similarly, a warrantor of a statement $E$, does not need any name. You
can call it ``the warrantor of $E$'' (assuming that it exists, of course). So, what matters is not ``who is a warrantor of 
$E$~?'', but only ``is there a warrantor of $E$~?''. 

This is of course why assumptions are most often anonymous, because, as explained above, an assumption is nothing other than the
declaration of a warrantor. If we have a named hypothesis such as ``Assume $E$ ($\zeta$).'' in a text, it is of course very
handy to refer to it by its name $\zeta$ (that in general is just a number), 
but it is by no means required, because in the scope of this assumption, $\zeta$ is the 
name of the unique warrantor of $E$, who does not need any name because it is unique. 
There is even no need to refer to it in general. Remark that if we had the funny idea
of assuming the same statement $E$ twice with distinct names $\zeta$ and $\eta$, it would be quite ridiculous to refer preferably to
$\zeta$ instead of $\eta$. Do you think that one of these two hypotheses could be better than the other one for any
(strictly mathematical) purpose~?

Similarly, a theorem can be anonymous (even if it is usual to give it the name of a
person who found a proof of it for the first time). When a mathematician uses a theorem proved in another article or book, he does not need to
specify that in order to be suitable for the purpose of his own proof, it must have been proved in one way rather than
in another. The only thing that matters is that the theorem was {\em reliably} proved at least once. Also remark that if a
mathematician finds a better proof of an already proved statement (better because shorter, more conceptual, using better
    tools, and so on\dots) the theorem keeps the name of the original mathematician. In a
sense, the name is not attributed with respect to the quality of the proof given, but only according to anteriority, that
is to the insurance that the warrantor exists. Formally, the name is given to the warrantor, not to the proof. 

It is a fact that proofs can be written in a lazy way, or more precisely that we do not write the proof of a theorem in
the same way depending on the audience. If you write a proof for undergraduate students, you will put more details than
if you write it for researchers. This just means of course that researchers have a better ability to find out the
missing (or implicit) parts of the proof. In both cases, it does not matter how people fill theses gaps, provided that they fill
them. This is again just because of the uniqueness of the warrantor. 

Furthermore, proofs are so lazily written in general that some necessary (formal) constructs are always omitted. In
other words, these constructs have no equivalent in the vernacular. As an example, assume that $p$ is a proof of an
implication $E\implies F$ (in any context), and that $q$ is a proof of $E$ (in the same context). Formally, a proof of $F$
is then the applicative term $p(q)$ (or something similar, in any case, a kind of ``modus ponens operator'' applied to 
    $p$ and $q$). 
There is nothing in the vernacular that corresponds to this
$p(q)$. We content ourself with the fact that $E$ and $E\implies F$ having become ``obvious'' at some point of the discourse,
$F$ is also obvious at this point. This is of course possible because we do not have to be precise on \emph{how} we prove this
fact, which is again a consequence of the uniqueness of the warrantor.
There is a mental mecanism corresponding to the possibility of ``applying'' $p$ to $q$, 
but nothing to writing it down. There are several other example of this phenomenon, and almost half of the formal proof
language constructs (that we do not need to describe here) is concerned. 

The above arguments are actually \emph{symptoms} of the fact that proof-irrelevance is at work in mathematics.
They are not an \emph{explanation} of the necessity of this principle. We do not pretend that the argument we are now 
giving is \emph{the} explanation, but it is nevertheless the most explanatory we have at hand, despite the fact that we
are conscious that our argument is still disputable.

First of all, it is necessary that the reader accepts that the vernacular language of mathematics is a strongly typed
language. This seems to be in contradiction with Zermelo-Fraenkel set theory, but actually there is no contradiction
because Zermelo-Fraenkel set theory is not a formalization of the language of mathematics, but only some way of
\emph{coding} mathematical objects into the sole notion of set. The simple fact that a declaration such as ``\emph{Let
$X$.}'' is unacceptable (because the reader will for sure ask \emph{What is $X$~?}), shows that mathematical objects have
a type and are not only \emph{nude sets} as Zermelo-Fraenkel set theory could let us think.(\footnote{Remark that we also
    often write ``\emph{Let $X$ be a set.}'', which show that sets are only particular objects among many others.}) 
Furthermore, almost all
mathematical symbols are \emph{polysemic} (in other words, they have several different meanings), and the only tool we
have for resolving the ambiguities created by polysemy is to use the types of objects. Notice that these types are mainly
generated by the fact that we give a name to each new concept we define. 

That said, a picture of mathematics can be that of a four levels system, where the levels contain respectively types, 
(mathematical) objects, statements and warrantors. The Curry-Howard correspondance just says that
the pair (types,objects) has essentially the same behaviour as the pair (statements,warrantors), and the purpose of this
discussion is to explain why statements have \emph{at most} one warrantor, unlike types that can have several distinct elements. 

As is well known in topos theory,  any statement can be rewritten in the form of an equality whose left and right 
hand members contain no statement at all. As a consequence, the unique purpose of proofs is to prove equalities between
objects of the same type.(\footnote{An equality between two objects of distinct types (such as a matrix and a polynomial) 
    is meaningless.}) Now, if it is the case that we need to make a distinction between two warrantors of the same
statement, we would need a tool for ensuring equality (or non equality) between warrantors. This tool could not be anything other than
another notion of proof, because warrantors would be at least as complicated and undecidable as mathematical objects. 
But that way, we are clearly entering into an infinite regression, which is unacceptable. In other words,
it is our opinion that proof-irrelevance is a mandatory foundational principle.

Summarizing, warrantors are generally omitted in notations, and the context provides the insurance that they exist,
which is enough. As a consequence, the notation $\inf(A)$ for the greatest lower bound of $A$ makes no reference to a
warrantor (or proof) of the fact that $A$ is nonempty, so that the name $\zeta$ of the required warrantor does not
appear in the notation. Since it does not appear in the notation $\inf(A)$, it is also useless to declare it in front 
of $A\not=\emptyset$. This
is why, dependent conjunction and implication are written $E\et F$ and $E\implies F$ in the vernacular 
and not $(\zpE)\et F$
and $(\zpE)\implies F$. But they remain nevertheless dependent.

{\footnotesize
\begin{tabular}{lll}
M. Herrmann~:&1.& Institut Mathématique de Jussieu (IMJ), Projet Logique,\\
&&Université Paris-Diderot. France\\
&2.& BioEmergences, INAF - CNRS, ISC - Paris Ile de France, Avenue de la\\ 
&&Terrasse, 91190 Gif Sur Yvette, France -- funded on France BioImaging\\ 
&&Agence Nationale de la Recherche ANR-10-INBS-04-05\\
&3.& ISC-PIF, Paris, France\\
A. Prouté~:&1.& Institut Mathématique de Jussieu (IMJ), Projet Logique,\\
&&Université Paris-Diderot. France
\end{tabular}}

\begin{verbatim}
   matthieu.herrmann@math.univ-paris-diderot.fr
   alp@math.univ-paris-diderot.fr
\end{verbatim}

\end{document}